\documentclass[12pt]{amsart}
\input epsf
\makeatletter
\let\atopwithdelims\@@atopwithdelims
\let\over\@@over
\newtheorem{theorem}{Theorem}[section]

\newtheorem{corollary}[theorem]{Corollary}
\newtheorem{definition}[theorem]{Definition}

\newcounter{enumit}

\newcommand{\toto}{\longrightarrow}

\newcommand{\for}{\textrm{for}}

\newcommand{\sS}{{\mathcal S}}

\begin{document}
\title[A Simple Bijection for the Regions the Shi Arrangement]
{A Simple Bijection for the regions of the Shi Arrangement of Hyperplanes}
\author{Christos~A.~Athanasiadis}
\address{\hskip-\parindent Christos~A.~Athanasiadis\\
Mathematical Sciences Research Institute\\
1000 Centennial Drive\\
Berkeley, CA 94720}
\email{athana@msri.org}
\author{Svante~Linusson}
\address{\hskip-\parindent Svante~Linusson\\
Department of Mathematics \\ 
Stockholm University \\ 
S-106 91
Stockholm, SWEDEN}
\email{linusson@matematik.su.se}
\thanks{Both authors were supported by the Mathematical Sciences Research 
Institute, Berkeley, California. The second author was also supported by 
the Swedish Science Council. Research at MSRI 
is supported in part by NSF grant 
DMS-9022140.}
\begin{abstract}
The Shi arrangement ${\mathcal S}_n$ is the arrangement of affine hyperplanes in ${\mathbb R}^n$ of the form $x_i - x_j = 0$ or $1$, for $1 \leq i < j \leq n$. It dissects ${\mathbb R}^n$ into $(n+1)^{n-1}$ regions, as was first proved by Shi. We give a simple bijective proof of this result. Our bijection generalizes easily to any subarrangement of ${\mathcal S}_n$ containing the hyperplanes $x_i - x_j = 0$ and to the extended Shi arrangements. 
\end{abstract}

\maketitle

\section{Introduction}

A \textit{hyperplane arrangement} $\mathcal A$ is a finite set of affine 
hyperplanes in ${\mathbb R}^n$. The \textit{regions} of $\mathcal A$ are the 
connected components of the space obtained from ${\mathbb R}^n$ by removing the hyperplanes of $\mathcal A$. A classical example is provided by the \textit{braid arrangement} ${\mathcal A}_n$. It consists of the hyperplanes in ${\mathbb R}^n$ of the form $x_i = x_j$ for $1 \leq i < j \leq n$, i.e. the reflecting hyperplanes of the Coxeter group of type $A_{n-1}$. Its regions correspond to permutations of the set $[n] := \{1, 2,\ldots,n\}$.

\begin{figure}[htpb]
\center{\mbox{\epsfbox{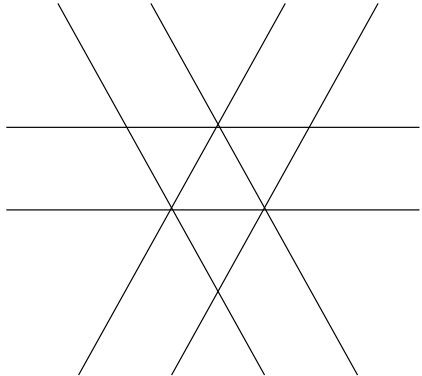}}} 
\caption{The Shi arrangement for $n=3$.}
\label{S3}
\end{figure}

\vspace{0.1 in}
A \textit{deformation} of ${\mathcal A}_n$ \cite{St1} is an arrangement each of whose hyperplanes is parallel to one of the hyperplanes of ${\mathcal A}_n$. We will be concerned with a deformation of ${\mathcal A}_n$ which has remarkable combinatorial properties. It is the \textit{Shi arrangement}, denoted by ${\mathcal S}_n$, and consists of the hyperplanes
\begin{equation}
\begin{tabular}{l}
$x_i - x_j = 0 \ \ \for \ \ 1 \leq i < j \leq n$\ and\\ 
$x_i - x_j = 1 \ \ \for \ \ 1 \leq i < j \leq n$ \,
\label{S_n}
\end{tabular} 
\end{equation}
in ${\mathbb R}^n$. Figure \ref{S3} shows ${\mathcal S}_3$ intersected with the plane $x_1 + x_2 + x_3 = 0$. Shi was the first to consider ${\mathcal S}_n$ in his investigation of the affine Weyl group of type $A_{n-1}$ \cite{Sh}. He used techniques from group theory to prove the following result.

\begin{theorem}
{\rm (Shi \cite[Cor.\ 7.3.10]{Sh})}
The number of regions of ${\mathcal S}_n$ is $(n+1)^{n-1}$.
\end{theorem}
The Shi arrangement was further studied enumeratively by Headley \cite{He1,He2}, Stanley \cite{St1,St2} and the first author \cite{Ath1,Ath2} and from the point of view of freeness in \cite{Ath3}. A simple proof of Theorem 1.1 was given in \cite[\S 3]{Ath1} \cite[\S 6.2]{Ath2} as an application of the ``finite field method'' of \cite{Ath1,Ath2}. It uses Zaslavsky's theory of counting regions \cite{Za}. Another simple proof of Theorem 1.1, using directly deletion and restriction, is implicit in \cite[Thm.\ 3.1]{Ath3}. See the first remark in the last section.

Our objective here is to give a \textit{simple bijective} proof. The only bijection already known is due to Pak and Stanley \cite[\S 5]{St1}. They established a correspondence between the regions of ${\mathcal S}_n$ and the \textit{parking functions} on $[n]$, which are well known to be counted by $(n+1)^{n-1}$. Although this correspondence is easy to define, a lot of effort is needed, as well as Shi's result itself, to prove that it is indeed a bijection (see the proof of \cite[Thm.\ 2.1]{St2}). Our bijection can also be stated in terms of parking functions but is different from that of Pak and Stanley. It generalizes easily to any arrangement between ${\mathcal A}_n$ and ${\mathcal S}_n$ as follows. Let $G$ be a simple graph on the vertex set $[n]$. We denote by ${\mathcal S}_{n,G}$ the arrangement
\begin{equation}
\begin{tabular}{l}
$x_i - x_j = 0 \ \ \for \ \ 1 \leq i < j \leq n,$\\ 
$x_i - x_j = 1 \ \ \for \ \ 1 \leq i < j \leq n, \ ij \in G$ 
\label{S_G}
\end{tabular} 
\end{equation}
in ${\mathbb R}^n$, first considered in \cite[\S 3]{Ath1} \cite[\S 6.2]{Ath2} and later in \cite{Ath3}. It specializes to ${\mathcal A}_n$ when $G$ is empty and to ${\mathcal S}_n$ when $G$ is the complete graph. Let ${\mathbb Z}_{n+1}$ denote the abelian group of integers modulo $n+1$ and let $H$ be the cyclic subgroup of ${\mathbb Z}_{n+1} ^n$ generated by $(1,1,\ldots,1)$. 
One can think of ${\mathbb Z}_{n+1} ^n$ as the set of all placements of 
$n$ distinct balls into $n+1$ identical boxes arranged cyclically.
\begin{theorem} 
The regions of ${\mathcal S}_{n,G}$ are in bijection with the cosets 
\begin{equation}
(a_1, a_2,\ldots,a_n) + H \in {\mathbb Z}_{n+1} ^n \, / H
\label{cos}
\end{equation}
which satisfy the following condition: given $i$, if $j$ is the smallest integer such that $i < j$ and $a_i = a_j$, then $ij \in G$.
\end{theorem}

Stanley \cite{St2} has generalized the correspondence of \cite[\S 5]{St1} to a bijection between the regions of the \textit{extended Shi arrangement}
\begin{equation}
x_i - x_j = -k+1, -k+2,\ldots,k \ \ \for \ \ 1 \leq i < j \leq n
\label{ext}
\end{equation}
and \textit{k-parking functions} on $[n]$, which are counted by $(kn+1)^{n-1}$. Our bijection also generalizes easily in this direction.

\vspace{0.1 in}
This paper is organized as follows: Section 2 contains our new proof of Theorem 1.1. In Section 3 we prove the more general Theorem 1.2 and derive some special cases, previously obtained with non-bijective methods. We also generalize our bijection to the extended Shi arrangements. In Section 4 we give explicitly the proof of Theorem 1.1 which follows from the methods of \cite{Ath3} and close with some open problems.

\section{The bijection}

We first describe our bijection in terms of parking functions. A \textit{parking function} on $[n]$ is a map $f: [n] \toto [n]$ such that for all $1 \leq j \leq n$, the cardinality of the set $f^{-1} ([j])$ is at least $j$. We also use the notation $f = (a_1, a_2,\ldots,a_n)$, where $a_i = f(i)$ for $1 \leq i \leq n$. Parking functions were first studied by Konheim and Weiss \cite{KW}. For the reason for the terminology ``parking function'' see \cite[\S 2.6]{Ha} \cite[\S 5]{St1}. An extensive literature is given in \cite{St2}. 

\smallskip
In order to describe the bijection we index the regions of $\sS_n$ as follows.
First, given a region $R$, we consider only the $x_i-x_j = 0$ hyperplanes and 
let $w = w_1w_2 \cdots w_n$ be the unique permutation of $[n]$ such that
$x_{w_1} > x_{w_2} > \cdots > x_{w_n}$ holds on $R$.  Second, we draw an arc 
$(i,j)$ in $w$ from $i$ to $j$ if $i<j$ and $x_i-x_j>1$ holds on $R$. Third, 
we remove any arcs ``containing'' another arc. In other words, if there is an 
arc $(j,k)$ then we remove any arcs $(i,l)$, $(i,k)$ or $(j,l)$ if 
$x_i > x_j > x_k > x_l$, $x_i > x_j$ or $x_k > x_l$ holds on $R$ respectively.
Clearly, these arcs are forced by the arc $(j,k)$ and hence redundant. The 
\textit{diagram} of $R$ is the resulting permutation of $[n]$ with arcs going 
rightwards from smaller to larger integers, with no arc containing another.

\begin{figure}[htpb]
\hspace{0 in} {\epsfysize=8mm\epsfbox{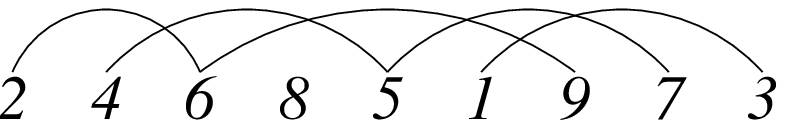}}
\caption{The diagram of a region of $\sS_9$}
\label{d9}
\end{figure}

\vspace{0.1 in}
\noindent
\textbf{Example.} The region of $\sS_9$ indexed by the diagram of Figure \ref{d9} is defined by the inequalities $x_2>x_4>x_6> \cdots >x_7>x_3$ and $x_2-x_4<1$, $x_2-x_6>1$, $x_4-x_6<1$, $x_2-x_8>1$ etc.

\smallskip
Note that for each diagram $\rho$, the arcs naturally determine a partition 
$\pi = \pi_{\rho}$ of $[n]$ into chains of increasing integers. In the 
example above, this partition is $\pi=269/457/8/13$. We say that the position 
of $m$ in $\rho$ is $j$ if $m = w_j$, i.e. if $m$ is the $j$th integer from 
the left which appears in $\rho$.  

\begin{definition} 
Let $\sigma_n$ be the map from the regions of $\sS_n$ to parking functions
on $[n]$ which sends the region with diagram $\rho$ to the function
\[f(i) = \mbox{the position in $\rho$ of the leftmost element in the chain 
containing $i$.}
\]
\end{definition}
This is clearly a parking function, so $\sigma_n$ is well defined. The region 
in our example is mapped by $\sigma_9$ to the 
parking function $(6,1,6,2,2,1,2,4,1)$.

\begin{theorem} [Main Theorem]
The map $\sigma_n$ is a bijection between the regions of $\sS_n$ and parking 
functions on $[n]$.  
\end{theorem}

\begin{proof}
We describe the inverse of $\sigma_n$ explicitly. Given a parking function $f$, we get the partition $\pi$ simply by placing $i$ and $j$ in the same 
block if $f(i) = f(j)$. The chains are obtained by listing the elements of 
each block in increasing order, from left to right. It remains to determine 
the permutation. To do so, we place the chains relative to each other one at a
time, in increasing order of their values under $f$. Assume that we have 
already placed the chains with values less than $j$ and are to 
place the chain with value $j$. Since $f$ is a parking function, there are at 
least $j-1$ elements already placed. We insert the leftmost element of the 
chain in position $j$, counting from the left. There is a unique way to place 
the other elements of the chain to the right without forming any pair of arcs 
with one containing the other. This braiding defines a diagram $\rho$ and 
hence a region $R$ of $\sS_n$. We leave it to the reader to check that this 
map is indeed the inverse of $\sigma_n$.
\end{proof}

\begin{figure}[htpb]
\hspace{0 in} {{\epsfysize=45mm\epsfbox{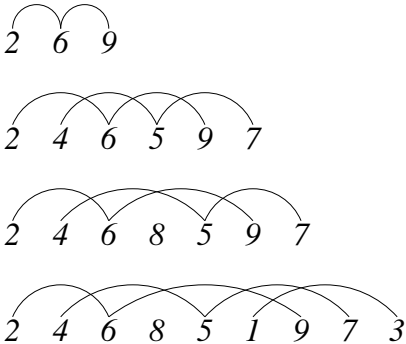}}}
\caption{Constructing the region $\sigma_9 ^{-1}(f)$}
\label{r9}
\end{figure}

Figure \ref{r9} illustrates the procedure to get back the region of $\sS_9$ 
from the parking function for our example.

\begin{figure}[htpb]
\center{\mbox{\epsfbox{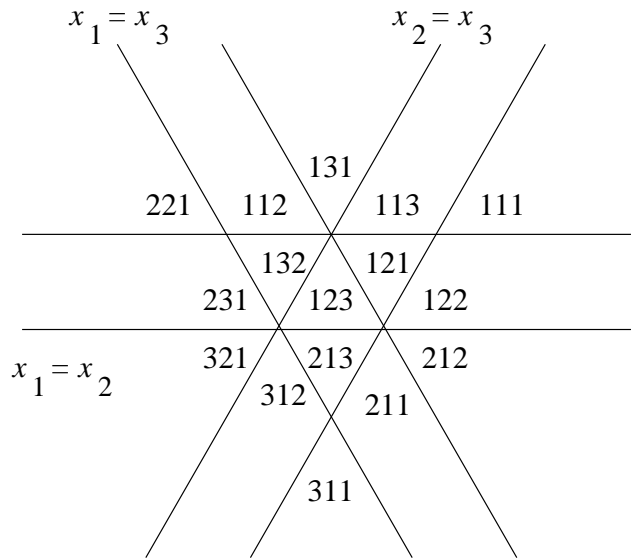}}} 
\caption{The bijection $\sigma_3$.}
\label{shi3}
\end{figure}

Figure \ref{shi3} shows the $16$ parking functions of length $3$ associated to the regions of $\sS_3$, according to $\sigma_3$. 

\vspace{0.1 in}
The fact that there are $(n+1)^{n-1}$ parking functions on $[n]$ follows from the observation, due to Pollack \cite[p. 13]{FR} and repeated by Haiman \cite[p. 28, 33]{Ha} and Stanley \cite[\S 2]{St3}, that every coset in ${\mathbb Z}_{n+1} ^n \, / H$ contains exactly one parking function. The following corollary proves Theorem 1.1.
\begin{corollary}
The map $\sigma_n$ induces a bijection between the regions of $\sS_n$ and elements of ${\mathbb Z}_{n+1} ^n \, / H$.
\end{corollary}

\section{Generalizations} 

In this section we generalize Theorem 2.2 to the arrangements between ${\mathcal A}_n$ and ${\mathcal S}_n$ and derive some special cases, previously obtained by other methods. Also, in a different direction, we give a generalization to the extended Shi arrangements.

\vspace{0.1 in}
\textit{Arrangements between ${\mathcal A}_n$ and ${\mathcal S}_n$.} 
Recall the definition of ${\mathcal S}_{n,G}$ given in (\ref{S_G}). A region $R$ of ${\mathcal S}_{n,G}$ can be represented as a permutation $w$ of $[n]$ together with a set of arcs, as in the case of ${\mathcal S}_n$. We now draw an arc $(i, j)$ in $w$ from $i$ to $j$ if $i < j$, $ij \in G$ and $x_i - x_j > 1$ holds on $R$. We remove all redundant arcs, as before, to get the \textit{diagram} of $R$. We define the map $\sigma_{n,G}$ from the regions of ${\mathcal S}_{n,G}$ to parking functions as in Definition 2.1.
\begin{theorem} 
The map $\sigma_{n,G}$ is a bijection between the regions of ${\mathcal S}_{n,G}$ and the parking functions $f = (a_1, a_2,\ldots,a_n)$ which satisfy the following condition: given $i$, if $j$ is the smallest integer such that $i < j$ and $a_i = a_j$, then $ij \in G$. 
\end{theorem}
\begin{proof}
Let $R$ be a region of ${\mathcal S}_{n,G}$ with diagram $\rho$. For a chain $i_1 < i_2 < \cdots < i_r$ in $\rho$ we have $i_{k-1} i_k \in G$ for all $1 < k \leq r$ by construction. Hence the associated parking function $f = (a_1, a_2,\ldots,a_n)$, for which $a_{i_1} = a_{i_2} = \cdots = a_{i_r}$, has the property stated in the theorem. The inverse of $\sigma_{n,G}$ is as in the special case of Theorem 2.2. 
\end{proof}

Theorem 1.2 follows immediately if we interpret parking functions as elements of ${\mathbb Z}_{n+1} ^n \, / H$, as in Section 2. As an application of Theorem 1.2 we obtain bijective proofs for two simple results from \cite{Ath1,Ath2}. The next theorem follows also from \cite[Cor.\ 3.6]{Ath3}.
\begin{theorem} 
{\rm (\cite[Thm.\ 3.4]{Ath1} \cite[Thm.\ 6.2.2]{Ath2})}
Suppose that the graph $G$ has the following property: if $1 \leq i < j < k \leq n$ and $ij \in G$ then $ik \in G$. Then the number of regions of ${\mathcal S}_{n,G}$ is the product
\begin{equation}
\prod_{1 < j \leq n} (n - d_j + 1),
\label{prod}
\end{equation}
where $d_j = \# \, \{i < j \ | \ ij \ \textrm{is not in} \ G \}$ for $1 < j \leq n$. 
\end{theorem}
\begin{proof}
Under the given assumption on $G$, the cosets (\ref{cos}) of Theorem 1.2 are exactly the ones that satisfy the following condition: if $i < j$ and $a_i = a_j$ then $ij \in G$. It suffices to show that the number of such cosets is the product (\ref{prod}). Indeed, fix a value for $a_1$ to break the cyclic symmetry and suppose we have chosen values $a_2,\ldots,a_{j-1}$ satisfying the condition. We want to choose $a_j \in {\mathbb Z}_{n+1}$ so that $a_j \neq a_i$ whenever $i < j$ and $ij$ not in $G$. These values $a_i$ are all distinct, since for two such $i_1 < i_2 < j$, $i_1 i_2$ is not in $G$ by the assumption on $G$ and hence $a_{i_1} \neq a_{i_2}$ by the choice of $a_{i_2}$. It follows that there are $d_j$ forbidden values for $a_j$ and hence $n - d_j + 1$ allowable ones. 
\end{proof}
\begin{theorem} 
{\rm (\cite[Thm.\ 5.6]{Ath1} \cite[Cor.\ 7.1.6]{Ath2})}
Let $G$ be the path $\{12, 23,\ldots,(n-1)n\}$. The number of regions of ${\mathcal S}_{n,G}$ is the sum
\[ \sum_{k=1}^{n} \, \frac{n!}{k!} \, {n-1 \choose k-1}. \]
\end{theorem}
\begin{proof}
Now a coset (\ref{cos}) satisfies the condition of Theorem 1.2 if and only if the entries which take any fixed value of ${\mathbb Z}_{n+1}$ form a string $a_i = a_{i+1} = \cdots = a_j$. There are ${n-1 \choose n-k} = {n-1 \choose k-1}$ ways to form $n-k+1$ such strings and $\frac{n!}{k!}$ ways to assign distinct values to them, modulo cyclic symmetry.    
\end{proof}

\vspace{0.1 in}
\textit{The extended Shi arrangements.} Following \cite{St2}, we denote by $\sS_n ^k$ the extended Shi arrangement (\ref{ext}). Stanley \cite{St2} has defined a \textit{k-parking function} on $[n]$ to be a sequence of positive integers $f = (a_1, a_2,\ldots,a_n)$ such that the unique increasing rearrangement $b_1 \leq b_2 \leq \cdots \leq b_n$ of the terms of $f$ satisfies $b_i \leq 1 + k(i-1)$ for all $i$. Thus a $1$-parking function is an ordinary parking function. He generalized the correspondence of \cite[\S 5]{St1} to a bijection between the regions of $\sS_n ^k$ and $k$-parking functions on $[n]$. He also noted that, in agreement to the $k=1$ case, $k$-parking functions on $[n]$ are in bijection with the cosets of the cyclic subgroup of ${\mathbb Z}_{kn+1} ^n$ generated by $(1, 1,\ldots,1)$, where ${\mathbb Z}_{kn+1}$ is the abelian group of integers modulo $kn+1$. Hence there are exactly $(kn+1)^{n-1}$ $k$-parking functions on $[n]$.

\vspace{0.1 in}
We now generalize the bijection $\sigma_n$ to treat the arrangements $\sS_n ^k$. We associate a diagram to a region $R$ of $\sS_n ^k$ as follows. First we consider only the hyperplanes $x_i - x_j = l$ for $-k+1 \leq l \leq k-1$ and let $y = y_1 \, y_2 \cdots y_{kn}$ be the unique permutation of the variables $x_i + m$, where $1 \leq i \leq n$ and $0 \leq m \leq k-1$, such that $y_1 > y_2 > \cdots > y_{kn}$ holds on $R$. We draw arcs in $y$ going rightwards from $x_i + m$ to $x_i + m-1$ for all $i$ and $m > 0$. Second, we draw an arc from $x_i$ to $x_j + k-1$ if $i < j$ and $x_i - x_j > k$ holds on $R$. Finally, we remove all arcs containing another arc and replace each variable $x_i + m$ by $i$. The arcs determine naturally a partition of the multiset $M_n ^k = \{1^k, 2^k,\ldots,n^k\}$ into chains of weakly increasing integers such that the elements of $M_n ^k$ equal to $i$ appear all in the same chain.

\begin{figure}[htpb]
\center{\mbox{\epsfbox{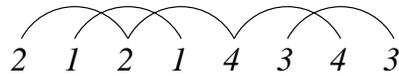}}}
\caption{The diagram of a region of $\sS_4 ^2$}
\label{d4-2}
\end{figure}

\vspace{0.1 in}
\noindent
\textbf{Example.} The diagram of Figure \ref{d4-2} represents the region of $\sS_4 ^2$ defined by the inequalities $x_2 + 1 > x_1 + 1 > x_2 > x_1 > x_4 + 1 > x_3 + 1 > x_4 > x_3$, $x_2 - x_4 > 2$ and $x_1 - x_3 < 2$. The corresponding partition of $M_4 ^2 = \{1, 1, 2, 2, 3, 3, 4, 4\}$ into chains is $2244/11/33$. 

\begin{definition}
We define the map $\sigma_n ^k$ by sending the region $R$ of $\sS_n ^k$ with diagram $\rho$ to the function $f = (a_1, a_2,\ldots,a_n)$ with
\[ a_i = \mbox{the position in $\rho$ of the leftmost element of the chain  
containing all $i$'s.}
\]
\end{definition}
As before, it is easy to check that $f$ is a $k$-parking function. For the region of Figure \ref{d4-2} we have $f = (2, 1, 6, 1)$.

\begin{theorem}
The map $\sigma_n ^k$ is a bijection between the regions of $\sS_n ^k$ and $k$-parking functions on $[n]$.  
\end{theorem}
\begin{proof}
We describe the inverse map of $\sigma_n ^k$, as in the proof of Theorem 2.2. Let $f = (a_1, a_2,\ldots,a_n)$ be a $k$-parking function and $b_1 \leq b_2 \leq \cdots \leq b_n$ be the unique increasing rearrangement of its terms. For each value $j$ of $f$, consider a chain $C_j$ of positive integers, listed from left to right in increasing order. The chain $C_j$ contains $k$ copies of $r$ if $a_r = j$ and none otherwise. Place the chains one at a time, in order of increasing value $j$. If $j = b_i > b_{i-1}$, then there are $k(i-1) \geq j-1$ elements listed before placing $C_j$, since $f$ is a $k$-parking function. Insert the leftmost element of $C_j$ in position $j$, counting from the left, and the other elements to the right so that no pair of arcs with one containing the other is formed. This defines the desired diagram, and hence region of $\sS_n ^k$. 
\end{proof}

\begin{figure}[htpb]
\hspace{0 in} {{\epsfysize=118mm\epsfbox{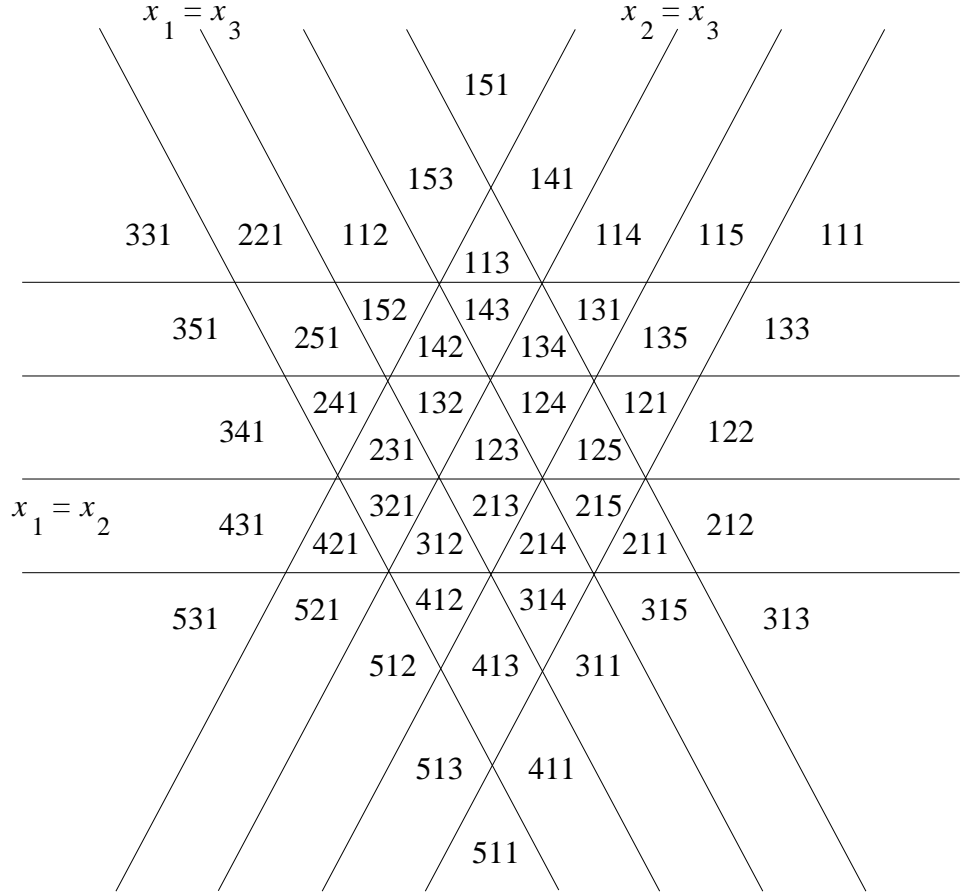}}} 
\caption{The bijection $\sigma_3 ^2$}
\label{shi2-3}
\end{figure}

Figure \ref{shi2-3} illustrates the bijection $\sigma_3 ^2$.

\section{Remarks and open problems}

1. The two other simple proofs of Theorem 1.1, mentioned in the introduction, come from computing the \textit{characteristic polynomial} \cite[\S 2.3]{OT} $\chi ({\mathcal S}_n, q)$ of ${\mathcal S}_n$. Assuming Shi's result, Headley \cite{He1,He2} was the first to show that $\chi ({\mathcal S}_n, q) = q (q-n)^{n-1}$. This formula is immediate \cite[Thm.\ 3.3]{Ath1} \cite[Thm.\ 6.2.1]{Ath2} once one uses the finite field method to interpret combinatorially the values of $\chi ({\mathcal S}_n, q)$ at large primes $q$. Zaslavsky's theorem \cite{Za} expresses the number of regions $r({\mathcal A})$ of a hyperplane arrangement ${\mathcal A}$ in ${\mathbb R}^n$ as $(-1)^n \chi ({\mathcal A}, -1)$ and yields Theorem 1.1. 

\smallskip
The computation of $\chi ({\mathcal S}_n, q)$ by deletion and restriction in \cite[Thm.\ 3.1]{Ath3} can be carried out on the level of the number of regions. It results in a naive inductive proof of Theorem 1.1, suitably generalized, which we describe explicitly next. Let ${\mathcal A}$ be a hyperplane arrangement and $H \in {\mathcal A}$ a distinguished hyperplane. The crucial and well known fact that we use below is that 
\begin{equation}
r({\mathcal A}) = r({\mathcal A}^{\prime}) + r({\mathcal A}^{\prime \prime}),
\label{d-r}
\end{equation}
where ${\mathcal A}^{\prime} = {\mathcal A} - \{H\}$ is the corresponding \textit{deleted arrangement} and ${\mathcal A}^{\prime \prime} = \{ H^{\prime} \cap H \ | \ H^{\prime} \in {\mathcal A}^{\prime} \}$ is the \textit{restricted arrangement} to $H$. 
Note that ${\mathcal A}^{\prime \prime}$ is an arrangement in 
the affine space $H$.
\begin{theorem}
For any integers $m \geq 0$ and $2 \leq k \leq n+1$, the arrangement
\begin{equation}
\begin{tabular}{l}
$x_1 - x_j = 0, 1,\ldots,m \ \ \for \ \ 2 \leq j < k,$\\
$x_1 - x_j = 0, 1,\ldots,m+1 \ \ \for \ \ k \leq j \leq n,$\\
$x_i - x_j = 0, 1 \ \ \for \ \ 2 \leq i < j \leq n$
\end{tabular}
\label{arr1}
\end{equation}
has $(n+m)^{k-2} (n+m+1)^{n-k+1}$ regions. In particular, for $m=0$ and $k=2$, ${\mathcal S}_n$ has $(n+1)^{n-1}$ regions.
\end{theorem}
\begin{proof}
We proceed by double induction on $n$ and $n-k$, the result being clear for $n=2$. The case $m=0$ and $k=n+1$ follows easily from the result for ${\mathcal S}_{n-1}$. Indeed, each of the $n^{n-2}$ regions of ${\mathcal S}_{n-1}$ in the space spanned by $x_2,\ldots,x_n$ determines a linear order of these variables and there are $n$ ways to form a region of (\ref{arr1}) by inserting $x_1$ in this order. We can now assume $2 \leq k \leq n$, since the arrangement (\ref{arr1}) having parameters $m \geq 1$ and $k = n+1$ coincides with (\ref{arr1}) having parameters $m-1$ and $k=2$. Consider the hyperplane $H$ of (\ref{arr1}) with equation $x_1 - x_k = m+1$. The corresponding deleted arrangement has the same form as (\ref{arr1}), with $k$ replaced by $k+1$, once one replaces $x_k$ by $x_1 - m - 1$. The restricted arrangement to $H$ has again the same form, with $n$ replaced by $n-1$ and $m$ replaced by $m+1$. The result follows by the induction hypothesis on these two arrangements and (\ref{d-r}).
\end{proof}

\vspace{0.1 in}
2. Let $k$ be any integer satisfying $1 \leq k \leq n$. The number of faces of ${\mathcal S}_n$ of dimension $k$ was shown \cite[Thm.\ 6.5]{Ath1} \cite[Cor.\ 8.2.1]{Ath2} to have the surprisingly simple combinatorial interpretation
\[ f_k({\mathcal S}_n) = {n \choose k} \, \# \, \{ f : [n-1] \rightarrow [n+1] \ \ | \ \ [n-k] \subseteq \textrm{Im} f \}, \]
where $\textrm{Im} f$ is the image of the map $f$. This formula reduces to Theorem 1.1 for $k=n$. The general case lacks a bijective proof and shows that the combinatorics of ${\mathcal S}_n$ is still not well understood. A similar interpretation was obtained for the extended Shi arrangements \cite[Thm.\ 8.2.2]{Ath2}.

\vspace{0.1 in}
3. Shi \cite{Sh2} has generalized Theorem 1.1 to the other irreducible crystallographic root systems. It would be interesting to find similar simple bijective proofs at least for the infinite families of type $B$, $C$ and $D$.

\newpage


\begin{thebibliography}{99}
%
\bibitem{Ath1}
\textsc{C.~A.~Athanasiadis}, 
Characteristic polynomials of subspace arrangements and finite fields, 
\textit{Advances in Math.} {\bf~122} (1996), 193--233.
%
\bibitem{Ath2}
\textsc{C.~A.~Athanasiadis}, 
\textit{Algebraic combinatorics of graph spectra, subspace arrangements and 
Tutte polynomials}, 
Ph.D. thesis, MIT, 1996.
%
\bibitem{Ath3}
\textsc{C.~A.~Athanasiadis},
On free deformations of the braid arrangement,
Preprint 526/1996, TU-Berlin, July 1996, 17 pages.
%
\bibitem{FR}
\textsc{D.~Foata and J.~Riordan}, 
Mappings of acyclic and parking functions,
\textit{Aequationes Math.} {\bf 10} (1974), 10--22.
%
\bibitem{Ha}
\textsc{M.~D.~Haiman}, 
Conjectures on the quotient ring by diagonal invariants,
\textit{J.\ Alg.\ Combin.} {\bf 3} (1994), 17--76.
%
\bibitem{He1}
\textsc{P.~Headley}, 
On reduced words in affine Weyl groups, 
in \textit{Proc.\ ``Formal Power Series and 
Algebraic Combinatorics (FPSAC) 1994''} 
(L.~J.~Billera, C.~Greene, R.~Simion, R.~Stanley, eds.),
DIMACS Series in Discrete Mathematics and Theoretical Computer Science,
Amer.\ Math.\ Soc. {\bf~24} (1996), pp. 225--232.
%
\bibitem{He2}
\textsc{P.~Headley}, 
\textit{Reduced expressions in infinite Coxeter groups},
Ph.D. thesis, University of Michigan, 1994.
%
\bibitem{KW}
\textsc{A.~G.~Konheim and B.~Weiss}, 
An occupancy discipline and applications,
\textit{SIAM J.\ Applied Math.} {\bf~14} (1966), 1266--1274.
%
\bibitem{OT}
\textsc{P.~Orlik and H.~Terao}, 
\textit{Arrangements of Hyperplanes},
Grundlehren 300, Springer-Verlag, New York, NY, 1992.
%
\bibitem{Sh}
\textsc{J.-Y.~Shi}, 
\textit{The Kazhdan-Lusztig cells in certain affine Weyl groups}, 
Lecture Notes in Mathematics, no.\ 1179, 
Springer-Verlag, Berlin/Heidelberg/New York, 1986.
%
\bibitem{Sh2}
\textsc{J.-Y.~Shi}, 
Sign types corresponding to an affine Weyl group,
\textit{J.\ London Math.\ Soc.} {\bf~35} (1987), 56--74. 
%
%
\bibitem{St1}
\textsc{R.~Stanley}, Hyperplane arrangements, interval orders and trees, 
\textit{Proc.\ Nat.\ Acad.\ Sci.} {\bf~93} (1996), 2620--2625.
%
\bibitem{St2}
\textsc{R.~Stanley}, Hyperplane arrangements, parking functions and tree
inversions, in \textit{Festschrift in Honor of Gian-Carlo Rota}, Birkh\"auser, 
Boston/Basel/Berlin, to appear.
%
\bibitem{St3}
\textsc{R.~Stanley}, Parking functions and noncrossing partitions, 
Preprint, August 1996.
%
\bibitem{Za}
\textsc{T.~Zaslavsky}, 
Facing up to arrangements: face-count formulas for partitions of space 
by hyperplanes,
\textit{Mem.\ Amer.\ Math.\ Soc.} vol.\ 1, no.~154, (1975).
%
\end{thebibliography}
\end{document}